\documentclass[12pt]{article}
\usepackage{amssymb, amsmath, amsthm}

\newcommand{\mr}[1]{\mathrm{#1}}

\newcommand{\z}{{\bf Z}}
\newcommand{\q}{{\bf Q}}
\newcommand{\zp}{{\bf Z}_p}
\newcommand{\qp}{{\bf Q}_p}

\newcommand{\cpn}{\z/p^n\z}

\newcommand{\zpn}{\zeta_{p^n}}

\newcommand{\mpn}{\mu_{p^n}}

\newcommand{\lpn}{\lambda_{n}}

\newcommand{\rt}[1]{\!\! \root p^n \of {#1}} 
\newcommand{\rti}[2]{\!\! \root p^{#1} \of {#2}}

\newcommand{\cs}{^{\times}} 
\newcommand{\cx}[1]{^{\times{#1}}}

\newcommand{\pp}[1]{{#1}\cs/{#1}\cx{p^n}}

\newcommand{\ra}{\rightarrow}

\newcommand{\lp}{\left(} 
\newcommand{\rp}{\right)}

\newcommand{\eq}[1]{(\ref{#1})}

\newcommand{\hatf}{D}

\newcommand{\ifs}{\mr{if \ }}
\newcommand{\tr}{Tr} 

\newcommand{\ord}{\Psi}

\newtheorem*{main}{Theorem}
\newtheorem{theorem}{Theorem}[section] 
\newtheorem{proposition}[theorem]{Proposition}
\newtheorem{lemma}[theorem]{Lemma} 
\newtheorem{corollary}[theorem]{Corollary}

\theoremstyle{definition}

\theoremstyle{remark}

\begin{document}
\title{Determination of Conductors from Galois Module Structure}
\author{Romyar T. Sharifi}
\date{October 2, 2001}
\maketitle
\begin{abstract}
    Let $E$ denote an unramified extension of $\qp$, and set $F = E(\zpn)$ 
    for an odd prime $p$ and $n \ge 1$.  We determine the conductors of the 
    Kummer extensions $F(\rt{a})$ of $F$ by those elements $a \in F\cs$ such 
    that $F(\rt{a})/E$ is Galois.  This follows from a comparison of 
    the Galois module structure of $F\cs$ with the unit filtration of $F$.
\end{abstract}
    
\section{Introduction}
\label{intro}
We consider ramification in those Galois extensions $K$ of an
unramified extension $E$ of $\qp$ which are also Kummer extensions of 
$F = E(\zpn)$ by the $p^n$th root of an element of $F\cs$, where $p$ 
is an odd prime and $\zpn$ is a primitive $p^n$th root of unity.   
We determine the 
conductors of the abelian subextensions of $F$ in $K$.  From these, 
the ramification groups of the entire extension $K/E$ are computable 
(see \cite{me-ram}, for example).

The two step solvable extensions we consider arise as fixed fields
of the kernels of certain upper triangular representations $\rho \colon 
G_{E} \ra GL_2(\cpn)$ of the absolute Galois group $G_{E}$ of $E$.
More specifically, for $\sigma \in G_{E}$ we should have
\[ 
  \rho(\sigma) = \lp \begin{array}{cc} 
    \chi^{s+t}(\sigma) & \kappa(\sigma) \\
    0 & \chi^{t}(\sigma)
  \end{array} \rp,
\]
where $\chi$ denotes the cyclotomic character 
with $\sigma(\zpn) = \zpn^{\chi(\sigma)}$ and $\kappa$ is a map which makes 
$\rho$ a 
homomorphism (in other words, $\kappa\chi^{-t}$ should be a $1$-cocycle of 
$G_{E}$ with values in $\cpn(s)$ for some $s$).
Such extensions can often be found as localizations of interesting
global extensions ramified only at primes above $p$ (see \cite{ihara}, 
for example).

We fix the fields $E$ and $F$ as above throughout the remainder of this
article.  Let $U_i$ denote the $i$th unit group of $F$ for $i \ge 0$.
We also let $(\,\cdot\,,\,\cdot\,)_{n,F}$ denote the $p^n$th norm residue 
symbol of $F$ (and similarly for $\qp(\zpn)$).  For $a \in F\cs$, we denote 
by 
$f_n(a)$ the conductor of $F(\rt{a})/F$ (considered as an integer) 
\cite[XV.2]{serre-lf}.  
Then $f = f_n(a)$ is the smallest nonnegative integer for which 
$(a,u)_{n,F} = 1$ for all $u \in U_f$.

For any integer $r$, we define a group $F^r$ by
\[
  F^r = \{ x \in F\cs \mid \sigma(x)x^{-\chi(\sigma)^r} \in F\cx{p^n}
  \mr{\ for\ all\ } \sigma \in G_E  \}.
\]
We shall compute the conductors $f_n(x)$ for all $x \in F^r$.  In particular,
in Section \ref{meta}, we shall prove the following result which gives the
answer for those $r$ with $r \not\equiv 0,1 \bmod p-1$.  

\begin{main}
  Let $r \not\equiv 0,1 \bmod p-1$, and let $t$ denote the smallest
  positive integer with $t \equiv 2-r \bmod p-1$.  Let $x \in
  F^r$ with $x \not\in F\cx{p}$.  Then we have that $f_m(x) =
  p^{m-1}t$ for any positive integer $m$ with $m \le n$.
\end{main}

The cases $r \equiv 0,1 \bmod p-1$ have more complicated statements
and are worked out in Sections \ref{metar0} and \ref{metar1}.
The case $r=0$, in which one considers extensions of $F$ by roots of elements 
of $E\cs$, was described using completely different methods in \cite{cm} 
for $E = \qp$ and \cite{me-cond} in general.

The determination of the conductors follows from the determination of 
generators of certain unit subgroups as Galois submodules of $F\cs$.  
More specifically, one can determine a simple formula for the norm
residue symbol for $(x,y)_{n,F}$ with $x \in F^r$ and $y \in U_1$
in terms of certain basic
symbols $(x,z)_{n,F}$ with $z$ running over a set of generators of $U_1$ 
as a Galois module.
Using this formula, one can evaluate $(x,y)_{n,F}$, with $y$ running
over the generators of the particular unit subgroups $U_{f-1}$ and
$U_f$, in order to determine that $f_n(x) = f$.

The determination of the generators can be done for all of the unit 
subgroups, but we do not follow this course here, as it is not 
needed for the determination of the conductors.

\vspace{.5ex}

\noindent{\em Acknowledgments.}\ A significant part of this article formed a 
part of the author's 
University of Chicago Ph.D. thesis \cite{me-thes}.  The author thanks 
Spencer Bloch and Dick Gross for their encouragement and advice.  He also 
thanks Hendrik Lenstra for a helpful early discussion.  The author was 
partially supported by NSF VIGRE grant 9977116 and an NSF Mathematical 
Sciences Postdoctoral Research Fellowship.

\section{Eigenspaces of the multiplicative group} \label{eigmult}

We maintain the notation of the introduction.
We may decompose $G = G_{F/\qp}$ with $p$ odd into a direct product of
cyclic subgroups 
\[ G = \Delta \times \Gamma \times \Phi, \] 
where $\Sigma = G_{F/E} = \Delta \times \Gamma$, the group $\Delta$ has order 
$p-1$, $\Gamma$ has order $p^{n-1}$ and $\Phi = G_{F/\qp(\zpn)}$.  

Let $\hatf$ denote the pro-$p$ completion of $F\cs$.  As the action 
of $\Delta$ on $\hatf$ is semisimple, $\hatf$ decomposes 
into a direct sum of eigenspaces for the powers of the cyclotomic character 
on $\Delta$.  For $r \in \z$, we let $D_r$ denote the eigenspace consisting 
of $x \in \hatf$ such that
\begin{equation*} 
  \delta(x) = x^{\omega(\delta)^r}, 
\end{equation*}
for all $\delta \in \Delta$,
where $\omega \colon \Delta \ra \z_p^*$ is the character with 
$\delta(\zeta_p) = \zeta_p^{\omega(\delta)}$.  
In other words, $D_r = \hatf^{\epsilon_r}$ where
\[ 
  \epsilon_r = % \frac{1}{p-1}
  \sum_{\delta \in \Delta} \omega(\delta)^{-r}\delta.
\]

It is the goal of this section to gain an understanding of the 
$A = \zp[\Gamma \times \Phi]$-module structure of $D_r$, viewed inside 
$\hatf$.
The result which follows shortly is clearly very important in this regard 
and results from an examination of the arguments of Greither 
in \cite{greither}.  
(See \cite{jannsen} for similar results of a less explicit nature.)

Let $\varphi \in \Phi$ denote the Frobenius element. 
Let $\gamma$ denote a generator of $\Gamma$ with $\chi(\gamma) \equiv 
1+p \bmod p^n$ (considering $\chi$ as a character of $G_{F/E}$), and let 
$\sigma$ be a generator of $\Sigma$ such that $\sigma^{p-1} = \gamma^{p-1}$. 
We let $\chi_{\sigma}$ denote a lift of $\chi(\sigma)$ to $\zp^*$.
Let $q$ denote the order of the residue field of $E$.
Let $\xi \in \mu_{q-1}(F)$ be such that $\tr_{\Phi} \xi = 1$
and the conjugates of $\xi$ form a normal basis of $F$ over 
$\qp(\zpn)$.

\begin{theorem} \label{struct}
   The $\zp[G]$-module $\hatf$ has a presentation
   \begin{equation} \label{multpres}
     \hatf = \langle u,v,w,\pi \mid \sigma v = v, \varphi \pi = \pi,
     N_{\Sigma} u = v^{1-\varphi}, N_{\Phi} u = \pi^{\sigma-1},
     w^{\chi_{\sigma}^{-1}\sigma -1} = u^{\varphi-1} \rangle,
   \end{equation}
   where $N_{\Sigma} \pi = p$, 
   $w \in U_1-U_2$, $u \in U_2-U_3$ and $v \equiv 1 + p\xi \bmod p^2$.
\end{theorem}

\begin{proof}
  In the notation of \cite{greither}, $u$ is the element 
  $u_{n-1}^{(1)}$, $w$ is $y_{n-1}^{r^{-1}}$, $v$ is $v_{n-1}^{(1)}$,
  and $\pi$ is $\pi_{n-1}$.  The properties of the elements
  listed at the end of the theorem
  follow from their definitions and the expansion of $\log u_{n-1}^{(1)}$
  found in \cite[p.\ 12]{greither}.
  
  In Theorem 3.3 of \cite{greither}, it
  is shown that there is an exact sequence
  \[ 1 \ra Q \ra \hatf \ra \mr{Ind}_{\Sigma}^{G} \mpn \ra 1, \]
  where $Q$ is the submodule of $\hatf$ generated by
  $u$, $v$ and $\pi$ subject to all but the last relation in 
  \eq{multpres}.  By their definitions, $u$ and $w$
  generate the kernel of the absolute norm on $F\cs$ as a 
  $\zp[G]$-module, and together with $v$, generate $U_1$.  Hence, the
  given elements generate $\hatf$.  Furthermore,
  Lemma 2.6 of \cite{greither} states that the final relation in 
  \eq{multpres} holds. 
  Therefore, $\hatf$ is a quotient of the group $P$ given by the
  presentation in \eq{multpres}.  It suffices to show that this 
  quotient becomes an isomorphism modulo $Q$.
  This follows from the following isomorphisms of 
  $\zp[G]$-modules:
  \[ 
    P/Q \cong \langle \bar{w} \mid \bar{w}^{\chi_{\sigma}^{-1}\sigma-1} = 
    1 \rangle \cong \mr{Ind}_{\Sigma}^G \mpn \cong \hatf/Q, 
  \]
  where the middle step is Lemma 2.4 of \cite{greither}.
\end{proof}

Note that, as a corollary, $N_{\Phi} w = \zpn$ for a proper choice
of $\zpn$, which we make.
Theorem \ref{struct} also gives us the structure of the groups $D_r$.

\begin{corollary} \label{eigenspaces}
    We maintain the notation of Theorem \ref{struct} and set $u_r = 
    u^{\epsilon_r}$, $\pi_0 = \pi^{\epsilon_0}$ and $w_1 = w^{\epsilon_1}$.
    If $r \not\equiv 0,1 \bmod p-1$, then $D_r$ is a free 
    $A$-module on $u_r$.  Furthermore, we
    have the $A$-module presentations
    \begin{equation} \label{D0}
      D_0 = \langle \pi_0,u_0,v \mid \gamma v = v, \varphi \pi_0 = 
      \pi_0, N_{\Gamma} u_0 = v^{1-\varphi}, N_{\Phi} u_0 = \pi_0^{\gamma-1} 
      \rangle
    \end{equation}
    and
    \begin{equation} \label{D1}
      D_1 = \langle u_1,w_1 \mid w_1^{(1+p)^{-1}\gamma-1} = 
      u_1^{\varphi-1} \rangle.
    \end{equation}
\end{corollary}
    
For $t \ge 1$, set $V_{t,r} = U_t \cap D_r$ and 
$V'_{t,r} = V_{t,r}-V_{t+1,r}$. 
We now briefly study the module structures of the $V_{t,r}$.

\begin{lemma} \label{eigmod}
    We have $V_{t,r}/V_{t+p-1,r} \cong {\bf F}_q$ for every $t \ge 0$, and
    $V'_{t,r} \neq 0$ if and only if $t \equiv r \bmod p-1$.        
\end{lemma}

\begin{proof}
    We sketch the proof.  Set $\lpn = 1-\zpn$, and
    let $\alpha = 1 + x\lpn^t$ for some $x \in U$.  For $\delta \in \Delta$,
    we have
    \[
       \delta(\alpha) % = 1 + (1-\zpn^{\omega(\delta)})^t\delta(x)
       \equiv 1+\omega(\delta)^t x \lpn^t \mod \lpn^{t+1}.
    \]
    Then
    \begin{equation} \label{alphaeps} 
      \alpha^{\epsilon_r} \equiv 1 - \sum_{\delta \in \Delta}
      \omega(\delta)^{t-r}x \lpn^t \mod \lpn^{t+1},
    \end{equation}
    and this is $1$ modulo $\lpn^{t+1}$ unless $t \equiv r \bmod p-1$.
    On the other hand, if $t \equiv r \bmod p-1$,
    then equation \eq{alphaeps} implies that
    \[
      \alpha^{\epsilon_r} \equiv \alpha \mod \lpn^{t+1}.
    \]
    The assertions of the lemma follow immediately.
\end{proof}

From now on, we set $V_t = V_{t,r}$ and $V'_t = V'_{t,r}$ if 
$t \equiv r \bmod p-1$.

\begin{lemma} \label{moving}
  Let $z \in V'_t$.  If $p \nmid t$, then
  $z^{\gamma-1} \in V'_{t+p-1}$.  Otherwise, $z^{\gamma-1} \in
  V_{t+2(p-1)}$.
\end{lemma}

\begin{proof}
  From the ramification groups of $F/E$ in the lower numbering
  \cite[IV.4]{serre-lf}, we see that $v(\gamma(\pi)/\pi - 1) = p-1$,
  where $v$ is the valuation on $F$ and $\pi$ is any prime element of
  $F$.  Hence, for any $x\in U_t$ with $p \nmid v(x-1)$, we have
  $v(\gamma(x)-x) = t+p-1$.  On the other hand, if $p \mid v(x-1)$
  then $v(\gamma(x)-x) > t+p-1$.  Since $z^{\gamma-1} \in D_t$, the
  result follows from Lemma \ref{eigmod}.
\end{proof}

\begin{lemma} \label{repeat}
  Let $z \in V'_t$.  Let $i$ denote the
  smallest nonnegative integer with $i \equiv t \bmod p$.
  For $j \ge 1$, we have that
  \begin{equation*}
    z^{{(\gamma-1)}^j} \in V_{t+(j+\left[\frac{j-i}{p-1}\right])(p-1)},
  \end{equation*}
\end{lemma}

\begin{proof} 
  Lemma \ref{moving} shows that $V_a^{\gamma-1} \subseteq
  V_{a+p-1}$ for all $a \equiv t \bmod p-1$.  Hence
  \begin{equation} \label{est1}
    z^{(\gamma-1)^j} \in V_{t+j(p-1)}. 
  \end{equation}
  Furthermore, Lemma \ref{moving} implies that 
  \[ z^{(\gamma-1)^i} \in V'_{t+(p-1)i} \] 
  and
  \[ y^{(\gamma-1)^{p-1}} \in V_{a+p(p-1)} \] 
  for $y \in V'_a$ with $a \equiv t \bmod p-1$.  Together these show that
  for every $p-1$ applications of $\gamma-1$ after the $i$th and
  starting with the $(i+1)$st, we can add an additional $p-1$ to the
  subscript in \eq{est1}, which yields the lemma. % \eq{exact}.
\end{proof}

The following easy lemma will also be useful.

\begin{lemma} \label{powers}
  Let $a$, $b$ and $j$ be nonnegative integers with $1 \le j \le n-1$
  and $a < p^{n-j}$.
  Then
  \[ p^jV_{a+b} \subseteq V_{p^ja+b}. \]
\end{lemma}

\section{The norm residue symbol on the eigenspaces}

We consider the norm residue pairing on the eigenspaces.
As $F^r \subseteq D_{r} F\cx{p^n}$, it will prove useful
to consider the restriction of the norm residue symbol to
elements lying in eigenspaces. 
Let $(\,\cdot\,,\,\cdot\,)$ denote the pairing induced by 
$(\,\cdot\,,\,\cdot\,)_{n,F}$ on $\hatf \times \hatf$.

\begin{lemma} \label{pair}
  The pairing
  \begin{equation*} 
    D_r \times D_s \ra \mpn 
  \end{equation*}
  induced by the Hilbert norm residue symbol is trivial unless $s
  \equiv 1-r \bmod p-1$ and nondegenerate otherwise.
\end{lemma}

\begin{proof}
  Let us evaluate $(a_r,a_s)$ with $a_r \in D_r$ and $a_s \in D_s$.
  Let $i = \omega(\delta)$ for a generator $\delta$ of $\Delta$.  Then
  we have
  \begin{equation*} 
    (a_r,a_s) = \delta(a_r,a_s)^{i^{-1}}. 
  \end{equation*}
  By Galois equivariance of the norm residue symbol, the last term
  equals
  \begin{equation*} 
    (\delta(a_r),\delta(a_s))^{i^{-1}} = 
    (a_r^{i^r},a_s^{i^s})^{i^{-1}}
    = (a_r,a_s)^{i^{r+s-1}}. 
  \end{equation*} 
  Hence either $i^{r+s-1} \equiv 1 \bmod p$ or $(a_r,a_s) = 1$.
  However, as $i$ has order dividing $p-1$, the former condition is
  exactly that $i^{r+s-1} \equiv 1 \bmod p^n$.  This proves the first
  statement, and the second follows by nondegeneracy of the norm
  residue symbol.
\end{proof}

From now on, we fix $r$ and let $s = 1-r$.
We have the following easy corollary of Lemmas \ref{eigmod} and \ref{pair}.

\begin{corollary} \label{condmod}
  Let $x \in D_r \cap F\cs$.  Then either 
  \[ f_n(x) \equiv s+1 \bmod p-1 \] 
  or $f_n(x) = 0$.
\end{corollary}

Define the symbol 
$[\,\cdot\,,\,\cdot\,]$ with values in $\cpn$ by
\[ 
  (\alpha,\beta) = \zeta^{[\alpha,\beta]}
\]
for $\alpha$, $\beta \in \hatf$.

\begin{lemma} \label{resform}
  Let $x \in F^r$ and $y \in D_s$.  Set
  \begin{equation*} 
    \alpha = y^{f(\gamma)} 
  \end{equation*}
  with $f \in \zp[X]$.  Then
  \begin{equation*} 
    [x,\alpha] \equiv f((1+p)^s)[x,y] \bmod p^n. 
  \end{equation*}
\end{lemma}      

\begin{proof}
  We need only show that $[x,\gamma(y)] = (1+p)^{s}[x,y]$.  We have
  \[
    (x,\gamma(y)) = \gamma(\gamma^{-1}(x),y) = (x^{(1+p)^{-r}},y)^{1+p}
     = (x,y)^{(1+p)^{s}}.
  \]
\end{proof}

Let $T = \gamma-1 \in \zp[\Gamma]$.

\begin{corollary} \label{trivsyms}
  Let $x \in F^r$, and let $m$ be a positive integer.
  Choose $y \in D_s$ and let $i$ be such that the order of $[x,y]$ is
  $p^{n-i}$.  Set $j = \min\{m+i,n\}$.  Then
  \begin{equation*} 
    [x,y^{p^kT^{m-k}}] \equiv 0 \bmod p^j
  \end{equation*}
  for $0 \le k \le m$.
\end{corollary}

\begin{proof}
  We calculate:
  \begin{equation*} 
    [x,y^{p^kT^{m-k}}] = p^k((1+p)^s-1)^{m-k}[x,y] 
    \equiv 0 \bmod p^j. 
  \end{equation*}
\end{proof}

For $\alpha,\beta \in \hatf$, let 
\[ 
  \ord(\alpha,\beta) = \min\{k  \mid p^k[\alpha,\varphi^i\beta] = 0
  \mr{\ for\ all\ } i \}.
\]

\begin{lemma} \label{nondegen}
  Let $x \in F^r$ with $x \notin F\cx{p}$.  
  If $S$ is a set of generators for $D_s$ as
  an $A$-module, then $\ord(x,y) = n$ for some $y \in S$.
\end{lemma}
    
\begin{proof}
  Lemma \ref{pair} implies that $(x,z)$ is a primitive $p^n$th root of unity
  for some $z \in D_s$.  We may write $z$ as a product of powers of 
  elements of the form $\varphi^i \gamma^j y$ with $y \in S$, 
  $i,j \in \z$.
  Hence $(x,\varphi^i \gamma^j y)$ is a primitive $p^n$th root
  of unity for some such $i$ and $j$.
  By Lemma \ref{resform},
  this cannot happen unless $(x,\varphi^i y)$ is a primitive
  $p^n$th root of unity for some $i$ as well.
\end{proof}

\section{Conductors for $r \not\equiv 0,1 \bmod p-1$} \label{meta}

We assume in this section that $r \not\equiv 0, 1 \bmod p-1$.  
The goal of this section is the determination of
the conductor of the extensions $F(\rt{x})/F$ for $x \in F^r$.  
By Theorem \ref{struct}, $D_s$ is a free $A$-module of 
rank 1.  We identify $D_s$ with $A$ via the isomorphism 
\[
  A \overset{\sim}{\ra} D_s,\ \alpha \mapsto u_s^{\alpha}
\]
of $A$-modules. 
Each submodule $V_i$ with $i \equiv s \bmod p-1$ 
becomes identified with an ideal of $A$ (independent of $u_s$) 
and by abuse of notation we treat them as equal. 
In addition, setting $\bar{F}^r = F^r/F\cx{p^n}$, we see that 
\[ \bar{F}^r \cong \cpn[\Phi], \]
as $D_r$ is also free of rank $1$.

Note that
\begin{equation*} 
  \zp[\Gamma] \cong \zp[Y]/((Y+1)^{p^{n-1}}-1). 
\end{equation*}
We shall require the following easy lemma.

\begin{lemma} \label{index}
  In $\zp[Y]$, the polynomial $(Y+1)^{p^{n-1}}-1$ is contained in the
  ideal
  \begin{equation*} 
    (p^{n-1}Y,p^{n-2}Y^p,\ldots,Y^{p^{n-1}}). 
  \end{equation*}
\end{lemma}

% \begin{proof}
%   It suffices to verify that if $p^m \le k < p^{m+1}$ with
%   $0 \le m < n-1$, then $v_p\lpb \binom{p^{n-1}}{k} \rpb \ge n-1-m$,
%   where $v_p$ denotes the $p$-adic valuation.  We note that, in fact,
%   \begin{equation*} 
%     v_p\lpb \binom{p^{n-1}}{k} \rpb= n-1-v_p(k) 
%   \end{equation*} 
%   for $1 \le k \le p^{n-1}$, which implies the claim.
% \end{proof}

Let $t$ denote the smallest positive integer congruent to $s+1$ modulo
$p-1$.
To determine the conductor $f_n(x)$ for $x \in F^r$, it suffices to look at
two ideals of $A$.

\begin{proposition} \label{2ideals}
  Let $m$ be a nonnegative integer with $m \le n-1$.
  The ideal $V_{p^{m}t-1}$ of $A$ contains an element $a$ of
  the form
  \begin{equation*} 
    a = p^{m} + \sum_{k=1}^{m} d_k p^{m-k} 
    T^{p^{k-1}t-1}
  \end{equation*}
  with $d_k \in \zp[\Phi]$ for $1 \le k \le m$.  
  Furthermore, we have 
  \[
    V_{p^{m}t+p-2} = (p^{m+1},p^{m}T, p^{m-k}T^{p^{k-1}t} 
    \mid 1 \le k \le m).
  \]
\end{proposition}

\begin{proof}
  We begin with the second statement.  
  By Lemma \ref{index}, the ideal
  \[ 
    I = (p^{m+1}, p^mY, p^{m-k}Y^{p^{k-1}t} \mid 1 \le k \le m)
  \]
  %$\zp[Y]$ generated by $p^{m+1}$, $p^{m}Y$ and
  %$p^{m-k}Y^{p^{k-1}t}$ for $1 \le k \le m$ 
  of $\zp[Y]$ contains $(Y+1)^{p^{n-1}}-1$.  
  Hence, the image $\bar{I}$ of $I$ in $\zp[\Gamma]$ has 
  %index equal to the index of $I$ in $\zp[Y]$, 
  \[ [\zp[\Gamma]:\bar{I}] = [\zp[Y]:I] = p^h, \]
  where
  \begin{multline*}
    h = m+1+(t-1)m+(p-1)t(m-1)+ \ldots + p^{m-2}(p-1)t \\=
    (p^{m-1}+\ldots+p+1)t+1.
  \end{multline*}
  Let $J = \bar{I}[\Phi]$ as an ideal of $A$.  
  Then Lemma \ref{eigmod} implies that 
  \[
    [V_{t-1}:V_{p^{m}t+p-2}] = q^h = [A:J].
  \]
  
  To prove the second statement, we are left only 
  to verify the claim that $J \subseteq V_{p^{m}t+p-2}$.
  Note that $1 \in V_{t-1}$, so clearly $p^{m+1} \in V_{p^mt+p-2}$.
  By Lemma \ref{moving}, we see that
  \[
    p^{m}T \in p^{m}V_{t+p-2}. 
  \]
  Similarly, Lemma \ref{repeat} yields
  \[
    p^{m-k}T^{p^{k-1}t} \in p^{m-k}V_{p^kt+p-2} 
  \]
  for $1 \le k \le m$.
  By Lemma \ref{powers}, we have
  \[
    p^{m-k}V_{p^kt+p-2} \subseteq V_{p^mt+p-2},    
  \]
  for $0 \le k \le m$, proving the claim.
      
  For the first statement, we proceed by induction, the case of $m = 0$
  being obvious.  Assume that we have proven the proposition for $m \le 
  n-2$.  Let $i = p^mt-1$.
  Then we have found $b \in V'_i$ of the desired 
  form.
  We claim that 
  \begin{equation} \label{somed}
      c = d T^i \in V'_{pi}
  \end{equation}
  for any
  $d \in \zp[\Phi]$ which is not a multiple of $p$.  By the second statement, 
  $c \notin V_{pi+2(p-1)}$, and by Lemma \ref{repeat} we have
  $c \in V_{pi}$.
  If it were that $c \in V'_{pi+p-1}$, then by Lemmas \ref{eigmod} and
  \ref{moving} we would have
  $d T^{i-1} \in V'_{pj}$
  with $j \le i-p+1$, but this would contradict 
  Lemma \ref{repeat}.
  Having demonstrated the claim, we see that 
  $a = pb + c \in V'_{pi+p-1}$ for an appropriate choice of 
  $d$ in \eq{somed} by Lemma \ref{eigmod}
  and the freeness of $D_r$ as a $\zp[\Phi]$-module.  
  It is clear that $a$ is the desired element for the case $m+1$ of the 
  first statement.
\end{proof}

We now determine the desired conductors.

\begin{theorem} \label{metcond}
  Let $r \not\equiv 0,1 \bmod p-1$, and let $t$ denote the smallest
  positive integer with $t \equiv 2-r \bmod p-1$.  Let $x \in
  F^r$ with $x \not\in F\cx{p}$.  Then we have that $f_i(x) =
  p^{i-1}t$ for any positive integer $i$ with $i \le n$.
\end{theorem}

\begin{proof}
  By Corollary \ref{trivsyms}, we have $[x,\varphi^jb] \equiv 0 \bmod p^i$ 
  for every $j$ and each
  generator $b$ of $V_{p^{i-1}t + p-2}$ listed in Proposition 
  \ref{2ideals} (with $m=i-1$).  
  On the other hand, by Lemma \ref{nondegen}, we may choose 
  $j$ such that for $u = \varphi^j(u_s)$, the symbol
  $(x,u)$ is a primitive $p^n$th root of unity.  In this case, we 
  have
  \begin{equation*} 
    [x,u^a] = [x,u^{p^{i-1}}] = p^{i-1}[x,u] \not\equiv 0 \bmod p^i, 
  \end{equation*}
  where $a \in \zp[\Gamma]$ is the element of Proposition \ref{2ideals}.
  Noting Corollary \ref{condmod}, we see that
  \begin{equation*} 
    f_i(x) = f_n(x^{p^{n-i}}) = p^{i-1}t.
  \end{equation*}
\end{proof}

\section{Conductors for $r \equiv 0 \bmod p-1$} \label{metar0}

Over the final two sections, we determine the conductors for 
the remaining values of $r$.  For an element $x \in F^r$, we will
first determine $f_n(x)$ in terms of the values 
$\ord(x,y)$, where $y$ runs over a set of $A$-module generators of
$D_{s}$.  Then, we shall explicitly describe these values in terms 
of $x$.  The situation is complicated by the existence of
multiple generators of $\bar{F}^r$ and $D_{s}$.  We shall twice appeal to 
the case of $r=0$ from \cite{me-cond}. 

In this section, we let $r \equiv 0 \bmod p-1$.
Recall the notation of Corollary \ref{eigenspaces}. 

\begin{lemma} \label{newgen}
  For $n \ge 2$, there exists an element $c \in \zp[\Phi]$ such that $y_1 = 
  u_1w_1^{pc}$ lies in $V'_{2p-1}$ and has image generating
  $V_{2p-1}/V_{3p-2}$ as a $\zp[\Phi]$-module.
\end{lemma}

\begin{proof}
  Recall from Theorem \ref{struct} that $w_1 \in V_1'$ and $u_1 \in V_p$.
  From this, Lemma \ref{moving}, Corollary \ref{eigenspaces} and the relation
  \[ w_1^{\gamma-1} = w_1^p u_1^{(\varphi-1)(1+p)} \]
  found in \eq{D1}, we see that $V_p/V_{3p-2}$ 
  is a $\z/p\z[\Phi]$-module generated by (the images of) $w_1^p$ and $u_1$.
  Now, the submodule $M$ generated by $w_1^{\gamma-1}$ and $w_1^p$ does not
  contain the submodule $N = V_{2p-1}/V_{3p-2}$, since the fact that 
  $N_{\Phi} w_1 = \zpn$ implies $N \cap N_{\Phi} M = 0$.  The conclusion now 
  follows.
\end{proof}

We remark that, with a little more effort, one can show that 
$c = \varphi^{-1} + aN_{\Phi}$ for some $a \in \zp$.
Clearly, $w_1$ and the element $y_1$ of Lemma \ref{newgen} generate $D_1$ as a 
$A$-module for $n \ge 2$.
Let us determine the generators of the relevant ideals of $D_1$.

\begin{proposition} \label{2idealsr1}
   Let $m$ be a positive integer with $m \le n-1$.
   \begin{enumerate}
      \item[a.] 
      The submodule $V_{p^{m}+p^{m-1}-1}$ of $D_1$ contains an element 
      of the form
      \[
        a_1 = p^{m}w_1 + \sum_{k=2}^{m} d_{k,1} 
        p^{m-k} T^{p^{k-1}+p^{k-2}-2} y_1
      \]
      with $d_{k,1} \in \zp[\Phi]$, and 
      \begin{equation} \label{submod} 
        V_{p^{m}+p^{m-1}+p-2} = (p^{m+1} w_1, p^{m-1}y_1, 
	p^{m-k}T^{p^{k-1}+p^{k-2}-1}y_1 \mid 2 \le k \le m).
      \end{equation}
      \item[b.]
      The submodule $V_{2p^{m}-1}$ of $D_1$ contains an element of
      the form
      \[
        a_2 = p^{m-1}y_1 + \sum_{k=2}^{m} d_{k,2} p^{m-k}
	T^{2p^{k-1}-2} y_1,
      \]
      with $d_{k,2} \in \zp[\Phi]$,
      and 
      \[
        V_{2p^{m}+p-2} = (p^{m+1} w_1,p^{m}y_1,p^{m-k}T^{2p^{k-1}-1}y_1 
        \mid 1 \le k \le m).
      \]	
   \end{enumerate}
\end{proposition}

\begin{proof}
  We proceed as in the proof of Proposition \ref{2ideals}.
  Let us focus on part a.
  Let $J$ be the submodule of $D_1$ given by the right hand side of 
  \eq{submod}.  
  This has index in $D_1$ equal to $q^{h}$, where
  \begin{multline*}
    h = m+1 + p(m-1) + (p^2-1)(m-2) + p(p^2-1)(m-3) + \ldots + 
    p^{m-3}(p^2-1) \\
    = 3+2p+2p^2+\ldots+2p^{m-2}+p^{m-1} = 
    \frac{p^m+p^{m-1}+p-3}{p-1}.
  \end{multline*}
  On the other hand, Lemma \ref{eigmod} implies that
  $[V_1:V_{i+p-1}] = q^{h}$ for $i = p^m+p^{m-1}-1$.
  Next, we remark that as $y_1 \in V'_{2p-1}$, we have
  \[
    p^{m-k}T^{p^{k-1}+p^{k-2}-1}y_1 \in p^{m-k}V_{p^k+p^{k-1}+p-2},
  \]
  for $2 \le k \le m$ by Lemma \ref{repeat}.
  Applying Lemma \ref{powers}, we have that $J \subseteq 
  V_{i+p-1}$ and hence equality by equality of the
  indices.  This proves the second statement of part 1.  
  
  For the first statement, the base case of $m = 1$ is obvious.  Assume that 
  we have proven the proposition for $m \le n-2$.  
  Then we have found $b \in V'_i$ %, where $i = p^m+p^{m-1}-1$, 
  of the desired form.  
  Using Lemmas \ref{moving}, \ref{repeat} and the second statement of
  the proposition, one can check (as in Proposition \ref{2ideals}) 
  that $c = d T^{i-1}y_1 \in V'_{pi}$ with 
  $d \in \zp[\Phi] - p\zp[\Phi]$.
  We conclude that $a_1 = pb + c \in V'_{pi+p-1}$ 
  is the desired element for some choice of $d = d_{m+1,1}$ .
  
  Part b follows from the same argument with the obvious modifications.
\end{proof}

If $n = 1$, we set $y_1 = u_1$.

\begin{proposition} \label{galcondr0}
  Let $x \in F^r-F\cx{p^n}$.  Let $i = \ord(x,y_1)$ and 
  $j = \ord(x,w_1)$.  Then
  \[
     f_n(x) = 
     \begin{cases}
	 p^n+p^{n-1} & \ifs i = n, \\
	 2p^i & \ifs 1 \le i \le n-1 \mr{\ and \ } j \le i+1, \\
         p^{j-1}+p^{j-2} & \ifs 2 \le j \le n \mr{\ and \ } i+2 \le j.
     \end{cases}
  \]
\end{proposition}

\begin{proof}
  For $n \ge 2$, with the exception of $j=1$ and $i=0$, this follows 
  directly from Proposition \ref{2idealsr1} after noting that 
  \[ V_{p^n+p^{n-1}-1} = pV_{2p^{n-1}-1} \]
  and similarly for $V_{p^n+p^{n-1}+p-2}$.  If $j = 1$ and $i = 0$,
  it follows from the fact that $w_1 \in V'_1$.  For $n=1$ and $i=1$, it 
  follows from $y_1 \in V'_p$ by Lemma \ref{eigmod}.
\end{proof}
  
Let $N_r \colon D_r \ra \bar{F}^r$ be defined by
\[ 
  N_r(x) = \prod_{i=1}^{p^{n-1}} \gamma^i(\bar{x})^{(1+p)^{-ir}}
\]
for $x \in D_r$ and $\bar{x}$ its image in $\pp{F}$.
By abuse of notation, we will use the same letter to denote an 
element of $\bar{F}^r$ and a chosen lift of it to $F^r$ (or conversely,
an element of $F^r$ and its image in $\bar{F}^r$).

Let $k \le n$ be maximal such that $r \equiv 0 \bmod p^{k-1}$.
In the notation of Theorem \ref{struct}, we set $t_r = N_r \pi_0$, 
$x_r = N_r u_0$ and $v_r = v^{p^{n-k}}$.  
Then $\bar{F}^r$ has  
presentation as a (multiplicative) $(\cpn)[\Phi]$-module given as
\begin{equation} \label{F0}
  \bar{F}^r = 
  \langle t_r,x_r,v_r \mid 
  t_r^{\varphi-1} = 1, N_{\Phi} x_r = t_r^{(1+p)^r-1}, 
  x_r^{p^{n-k}} = v_r^{1-\varphi}, v_r^{p^k} = 1 \rangle.
\end{equation}
Note that $k$ is maximal among $k \le n$ such that $p^k$ divides 
$(1+p)^r-1$.

For $\alpha \in \zp[\Phi]$, let $\nu(\alpha)$ denote the largest 
integer such that $\alpha \in p^{\nu(\alpha)}\zp[\Phi]$ (possibly 
infinite).  Note that $\nu$ is exactly the $p$-adic valuation on 
elements of $\zp$.

\begin{lemma} \label{symbolsr0}
  Let $\alpha \in \zp[\Phi]$ with $\nu(\alpha) = 0$.
  Let $\delta \in \zp[\Phi]$ with $\nu(\delta) = 0$ and
  $(\varphi-1)\delta \neq 0$.
  The following statements hold:
  \begin{align*}      
    &a.\  \ord(t_r,y_1) = n 
    &&b.\  \ord(t_r,w_1) = 0 \\
    &c.\  \ord(v^{\alpha},y_1) = n-1 
    &&d.\  \ord(v^{\alpha},w_1) = n \\
    &e.\  \ord(x_r^{\delta},y_1) = n-1 
    &&f.\  \ord(x_r^{\delta},w_1) = n
  \end{align*}
\end{lemma}

\begin{proof}
  We proceed case by case.
  \begin{list}{\emph{\alph{enumi}.}}{\usecounter{enumi}}
    \item[\emph{a.}] 
    This follows easily from Lemma \ref{nondegen} and
    part b (to be proven). 
    \item[\emph{b.}]   
    Set $\lambda = \sigma^i(\pi)$ for any $i$.  Since $N_{\Sigma} 
    \lambda = p$ from Theorem \ref{struct}, we have 
    $\lambda = (1-\zpn)\eta$ for some $\eta \in U_1$ with 
    $N_{\Sigma} \eta = 1$.  Then
    \[
      (\lambda,w_1)_{n,F} = (\lambda,\zpn)_{n,\qp(\zpn)} = 
      (\eta,\zpn)_{n,\qp(\zpn)} = 1, 
    \]
    the first step following from $\lambda^{\varphi-1} = 1$ and $N_{\Phi} 
    w_1 = \zpn$ and the last since $\zpn$ pairs trivially with any
    element in the kernel of the norm.
    The result now follows, as $t_r = N_r N_{\Delta} \pi$.
    \item[\emph{c.}] 
    In \cite[Theorem 8]{me-cond}, it was shown that
    $f_n(v^{\alpha}) = 2p^{n-1}$, and hence Proposition \ref{galcondr0} 
    forces that $\ord(v^{\alpha},y_1) = n-1$.
    \item[\emph{d.}] 
    This follows from part c and Lemma \ref{nondegen}.
    \item[\emph{e.}]
    By \eq{F0} and part c, we see that
    \begin{equation} \label{mostofxy} 
       (x_r^{\delta},\varphi^i y_1)^{p^{n-k}} = 
       (v^{\delta(1-\varphi)},\varphi^iy_1)^{p^{n-k}} 
    \end{equation}
    is a primitive $p^{k-1}$st root of unity for some $i$.
    Hence, we are done if $k \ge 2$.  If $k = 1$, then
    \[ 
      (x_r,N_{\Phi}y_1) = (N_{\Phi}x_r,y_1) = (t_r^{(1+p)^r-1},y_1)
    \]
    is a primitive $p^{n-1}$st root of unity by part a.
    Since $\ord(x_r,y_1) \le n-1$ by \eq{mostofxy}, we have equality.
    \item[\emph{f.}] 
    Since
    $x_r^{\delta}v^{(\varphi-1)\delta} \in F\cx{p^k}$ by \eq{F0}, this 
    follows from part d.
  \end{list}
\end{proof}  

The following describes the conductors of all elements of $F^r$
for $r \equiv 0 \bmod p-1$.

\begin{theorem} \label{condr0}
  Let $x = t_r^{\alpha} v^{\beta} x_r^{\delta}$ with
  $\alpha \in \zp$, $\beta \in p^{n-k}\zp[\Phi]$ and
  $\delta \in \zp[\Phi]$ satisfying either $\delta = 0$, or 
  $\nu(\delta) \le n-k-1$ and $(\varphi-1)\delta \neq 0$.
  Let
  \[
    i = n - \min\{ \nu(\alpha),\,\nu(\beta)+1,\,\nu(\delta)+1 \}
  \]
  and
  \[
    j = n - \max\{ \nu(p\beta-N_{\Phi}\alpha),0 \}.
  \]
  Then
  \[
    f_n(x) =
    \begin{cases}
      p^{i-1}(p+1) & \ifs i = 0,\ \mr{or\ }
        \delta = 0\ \mr{and\ } j < i,\\
      2p^i & \ifs i \le n-1\ \mr{and\ either\ }
        \delta \neq 0\ \mr{or\ } j = i,\\
      0 & \mr{otherwise}.
    \end{cases}	
  \]
\end{theorem}

\begin{proof}
  Part a of Lemma \ref{symbolsr0} yields that
  \[
    f_n(t_r^{\alpha}) = 
    \begin{cases} 
      p^{n-1}(p+1) & \ifs \nu(\alpha) = 0, \\
      2p^{n-\nu(\alpha)} & \ifs \nu(\alpha) \le n-1, \\
    \end{cases}
  \]
  for $\alpha \in \zp$.
  Parts c and d yield
  \[
    f_n(v^{\beta}) = 2p^{n-1-\nu(\beta)}
  \]
  if $\beta \in \zp[\Phi]$ satisfies $\nu(\beta) \le n-1$.
  Parts e and f yield
  \[
    f_n(x_r^{\delta}) = 2p^{n-1-\nu(\delta)}
  \]
  with $\delta$ satisfying $(\varphi-1)\delta \neq 0$ and $\nu(\delta) 
  \le n-1$.
  
  By the conditions on $\delta$ and $\beta$ that $\nu(\delta) \le 
  n-k-1$, 
  unless $\delta = 0$ and $\nu(\beta) \ge n-k$, we have
  \[
    f_n(v^{\beta}x_r^{\delta}) = 
    \max\{f_n(v^{\beta}),f_n(x_r^{\delta})\}.
  \]
  While $\ord(t_r,y_1^{\varphi-1}) = 0$, on the other hand, we have
  \[ 
    \ord(v^{\beta},y_1^{\varphi-1}) = n-1-\nu(\beta)
  \]
  by part c, unless $(\varphi-1)\beta = 0$, and
  \[
    \ord(x_r^{\delta},y_1^{\varphi-1}) = n-1-\nu(\delta) 
  \]
  by part e.
  Therefore, by Proposition \ref{galcondr0}, we have
  \[
    f_n(t_r^{\alpha}v^{\beta}x_r^{\delta}) =
    \max\{f_n(t_r^{\alpha}),f_n(v^{\beta}),f_n(x_r^{\delta})\}
  \]
  unless $\delta = 0$, $(\varphi-1)\beta = 0$ and 
  $\nu(\eta) = \nu(\beta) + 1 \le n-1$. 
  
  In this ``exceptional'' case, there exists 
  $c \in \zp$, unique modulo $p$, such that
  \[ \ord(t_r^{c\alpha} v^{\beta}, y_1) \le n-2-\nu(\beta). \]
  By Theorem \ref{struct}, we have
  $N_{\Phi} v \equiv 1+p \bmod p^2$ and, clearly,
  $t_r p^{-1} \in F\cx{p}$.  Thus, we see from
  \cite[Theorem 8]{me-cond} (or \cite[Theorem 6.1]{cm}) that
  \[ f_n(t_r^pN_{\Phi}v) = p^{n-2}(p+1) < f_n(N_{\Phi}v). \]
  The case statement of the theorem now follows.
\end{proof}

The use of \cite[Theorem 8]{me-cond} in Lemma 
\ref{symbolsr0} (aside from the case $(\varphi-1)\alpha = 0$ in part 
c) could have been avoided, but this would not have alleviated the need 
for its use in Theorem \ref{condr0}.  We remark that the methods used in
\cite{me-cond} are quite elementary, relying upon basic 
properties of the norm residue symbol and the Artin-Hasse law for the symbol 
$(\zpn,\,\cdot\,)$.

\section{Conductors for $r \equiv 1 \bmod p-1$} \label{metar1}

Note that $u_0 \in V'_{p-1}$ and $v \in V'_{p^{n-1}(p-1)}$. 

\begin{proposition} \label{2idealsr0}
  Let $m$ denote a positive integer with $m \le n$.
  \begin{enumerate}
  \item[a.]
  The submodule $V_{p^m-1}$ of $D_0$ contains an element of the form
  \[
    a_1 = p^{m-1}u_0 + \sum_{k=2}^{m} d_{k,1} p^{m-k} T^{p^{k-1}-1}u_0
  \]
  with $d_{k,1} \in \zp[\Phi]$, and if $m \le n-1$, then
  \[
    V_{p^m+p-2} = (v, p^m u_0, p^{m-k} T^{p^{k-1}}u_0 \mid 1 \le k \le m).
  \]
  \item[b.]
  The submodule $V_{p^n-1}$ of $D_0$ contains an element of the form
  \[
    a_2 = v + \sum_{k=2}^{n} d_{k,2} p^{n-k} T^{p^{k-1}-1}u_0
  \]
  with $d_{k,2} \in \zp[\Phi]$, and
  \begin{equation} \label{Vgens}
    V_{p^n+p-2} = (pv, p^n u_0, p^{n-k}T^{p^{k-1}}u_0 \mid 1 \le k \le 
    n-1).
  \end{equation}
  \end{enumerate}
\end{proposition}

\begin{proof}
  The proof of part a
  is virtually identical to those of Propositions \ref{2ideals} and 
  \ref{2idealsr0}.  
  Note that the key to proving the first statement is the fact that
  \begin{equation} \label{anote}
      dT^{p^{k-1}-1}u_0 \in V'_{p^k-p}
  \end{equation}
  for $2 \le k \le n$ and any $d \in \zp[\Phi]$ with $\nu(d) = 0$.
  
  We focus on part b.  Its second statement \eq{Vgens} follows the fact 
  that $pV_{p^{n-1}+p-2} = V_{p^n+p-2}.$ 
  As for the first statement, we begin by proving
  the claim that $N_{\Phi} a_1 \in V_{p^n+p-2}$, where $a_1$ is the 
  element of part a for $m=n$.  First, we remark that 
  \[ 
    N_{\Gamma} = \frac{(1+T)^{p^{n-1}}-1}{T} = \sum_{j=1}^{p^{n-1}} 
    \binom{p^{n-1}}{j} T^{j-1}. 
  \]
  Using \eq{Vgens} (and some simple congruences for
  binomial symbols), we see that
  \[ 
    N_{\Gamma} u_0 \equiv \sum_{k=1}^{n} p^{n-k} T^{p^{k-1}-1} u_0
    \bmod V_{p^n+p-2}.
  \]
  By \eq{D0}, we therefore have
  \[ 
    N_{\Phi} a_1 = N_{\Phi} a_1 - N_{\Gamma} N_{\Phi} u_0 
    \equiv \sum_{k=2}^{n-1} 
    N_{\Phi}(d_{k,1}-1)p^{n-k}T^{p^{k-1}-1} u_0 \bmod V_{p^n+p-2}.
  \]
  Assume that the claim does not hold, so that
  $N_{\Phi}(d_{k,1}-1) \not\equiv 0 \bmod p$ for some minimal $k$.
  Noting \eq{anote}, we see that $N_{\Phi} a_1 \in 
  V'_{p^n-p^{n-k+1}}$, contradicting $a_1 \in V'_{p^n-1}$.
  
  Consider the submodule
  \[  W = Au_0 \cap V_{p^n-1} = Aa_1 + V_{p^n+p-2} \]
  of $D_0$.
  From the claim, we see that $[W:V_{p^n+p-2}] = q/p$,
  whereas $[V_{p^n-1}:V_{p^n+p-2}] = q$ by Lemma \ref{eigmod}.
  Hence, by \eq{D0} and \eq{Vgens}, 
  there exists an element $b \in V'_{p^n-1}$ of the form
  \[
    b = v + \sum_{k=1}^{n} c_{k,2}p^{n-k}T^{p^{k-1}-1}u_0
  \]
  with $c_{k,2} \in \zp[\Phi]$.  We take $a_2 = b-c_{1,2}a_1$.
\end{proof}

Now let $r \equiv 1 \bmod p-1$.  The following is a direct consequence
of Proposition \ref{2idealsr0}.

\begin{proposition} \label{galcondr1}
  Let $x \in F^r$.  Let $i = \ord(x,u_0)$ and $j = \ord(x,v)$.  
  Then
  \[
     f_{n}(x) = 
     \begin{cases}
	 p^{n-1}(j(p-1)+1) & \ifs j \ge 1, \\
	 p^i & \ifs j=0,\ i \ge 1.
     \end{cases}
  \]
\end{proposition}

Let $k \le n$ be maximal such that $r \equiv 1 \bmod p^{k-1}$.
Let $x_r = N_r u_1$, $z_r = N_r w_1$ and $\kappa_r = 
\zeta_{p^n}^{p^{n-k}}$, and
recall that $N_{\Phi} w = \zpn$.  
From \eq{D1}, we obtain a presentation
\begin{equation} \label{F1}
  \bar{F}^r = \langle x_r,z_r,\kappa_r \mid 
  x_r^{\varphi-1} = z_r^{(1+p)^{r-1}-1}, N_{\Phi} z_r = \kappa_r^{p^{k-1}},
  \kappa_r^{p^k} = 1, \kappa_r^{\varphi-1} = 1 \rangle.
\end{equation}

We make the following useful remark.

\begin{lemma} \label{unnondegen}
  If $x \in F^r$ and $x \notin 
  F\cx{p}$, then the larger of $\ord(x,u_0)$ and $\ord(x,v)$ is at
  least $n-k$.
\end{lemma}

\begin{proof}
  We recall that for an element $x \in F\cs$ and $m \le n$, 
  the extension $F(\rti{m}{x})/F$ is 
  unramified if and only if $(x,u)_{m,F} = 1$ for all $u \in 
  U_1$.  Furthermore,
  $p^k$ is the degree of the maximal unramified extension of $D_r$ by 
  the $p^n$th root of an element of $F^r$, since such an extension
  must be abelian over $\qp$.  The result now follows.
\end{proof}

\begin{lemma} \label{symbolsr1}
  Let $\alpha \in \zp[\Phi]$ with $\nu(\alpha) = 0$.
  Let $\delta \in \zp[\Phi]$ with $\nu(\delta) = 0$ and
  $(\varphi-1)\delta \neq 0$.
  The following statements hold:
  \begin{align*}      
    &a.\  \ord(x_r^{\alpha},v) = 0
    &&b.\  \ord(x_r,u_0) = n-k \\
    &c.\  \ord(z_r^{\alpha},v) = 1 
    &&d.\  \ord(z_r^{\delta},u_0) = n \\
    &e.\  \ord(\zpn,v) = n
    &&f.\  \ord(\zpn,u_0) = 0
  \end{align*}
\end{lemma}

\begin{proof}
  Again, we proceed case by case.
  \begin{list}{\emph{\alph{enumi}.}}{\usecounter{enumi}}
    \item
    It suffices to consider $\alpha = 1$.
    We have (abusing notation)
    \begin{equation} \label{xa}
      [x_r,\varphi^i v] = [N_r u_1,\varphi^i v] = [N_r y_1,\varphi^i v] - 
      p[N_r w_1^c,\varphi^{i} v],
    \end{equation}
    with $c \in \zp[\Phi]$, by Lemma \ref{newgen}.
    Now apply Lemma \ref{resform} and the definition of $N_r$
    to see that the right hand side of \eq{xa} reduces to
    $p^{n-1}[y_1,a]$.
    The conclusion now follows from Lemma \ref{symbolsr0}c.
    \item 
    Note that $x_0 = v^{\varphi-1}$ by Corollary \ref{eigenspaces}.  
    Hence, for any $i$, we have
    \[
      (x_0,u_1^{\varphi^i}) = (v,u_1^{\varphi^{i-1}(1-\varphi)}) = 
      (v,w_1^{\varphi^{i-1}(1-(1+p)^{-1}\gamma)}) = 1
    \]
    by Lemma \ref{resform}.
    Since $x_{1-r}x_0^{-1} \in F\cx{p^k}$, we have that $\ord(x_{1-r}
    ,u_1) \le n-k$.  
    We conclude by remarking that
    \[ 
      (x_r,u_0^{\varphi^i}) = (u_1^{\varphi^{-i}},x_{1-r})
    \]
    and applying Lemma \ref{unnondegen}.
    \item 
    This follows from Lemma \ref{symbolsr0}d 
    in the same way that part a of this lemma followed from
    Lemma \ref{symbolsr0}c.
    \item 
    This follows from Lemma \ref{symbolsr0}f since
    \[ (z_r^{\delta},u_0) = (w_1,x_{1-r}^{\iota(\delta)}), \]
    where $\iota$ is the involution of $\zp[\Phi]$ defined by $\varphi 
    \mapsto \varphi^{-1}$.
    \item
    Recall that
    \begin{equation} \label{zetaform}
     [\zpn,x] = \frac{N_G x-1}{p^n}
    \end{equation}
    for any $x \in U_1$.
    Since $N_G v \equiv 1+p^n \bmod p^{n+1}$,
    we have the result. 
    \item
    This follows from \eq{zetaform} as well, since $N_G u_0 = 1$.
  \end{list}
\end{proof}  

The conductors of the elements of $F^r$ for $r \equiv 1 \bmod p-1$ are now 
described by the following.

\begin{theorem} \label{condr1}
  Let $x \in F^r$ with image $x_r^{\alpha}z_r^{\beta}\kappa_r^{\delta} 
  \in \bar{F}^r$, where $\alpha,\delta \in \zp$
  and $\beta \in \zp[\Phi]$ satisfies $(\varphi-1)\beta \neq 0$ if
  $\beta \neq 0$.  Let $i = k-\nu(\delta)$ and 
  $j = n-\min\{\nu(\beta),\nu(\alpha)+k\}$.
  Then
  \[
    f_n(x) =
    \begin{cases}
      	p^{n-1}(i(p-1) + 1) & \ifs i \ge 1, \\
	p^{j} & \ifs i \le 0 \mr{\ and\ } j \ge 1, \\  
	0 & \mr{otherwise}.
    \end{cases} 	
  \]
\end{theorem}

\begin{proof}
  Parts a and b of Lemma \ref{symbolsr1} yield 
  \[
    f_n(x_r^{\alpha}) = p^{n-k-\nu(\alpha)}
  \]
  if $\alpha \in \zp$ with $\nu(\alpha) \le n-k-1$ (and $0$ otherwise).
  Parts c and d yield 
  \[
    f_n(z_r^{\beta}) = p^{n-\nu(\beta)}
  \]
  if $\nu(\beta) \le n-1$ and $(\varphi-1)\beta \neq 0$.
  Part e yields 
  \[
    f_n(\kappa_r^{\delta}) = p^{n-1}((k-\nu(\delta))(p-1)+1)
  \]
  if $\delta \in \zp$ with $\nu(\delta) \le k-1$.
  
  Furthermore, we see from the values in Lemma \ref{symbolsr1}
  along with Proposition \ref{galcondr1} that
  \[
    f_n(x) = \max\{ f_n(x_r^{\alpha}),f_n(z_r^{\beta}),
    f_n(\kappa_r^{\delta}) \}
  \]
  with $\alpha$, $\beta$ and $\delta$ as in the theorem,
  unless perhaps if $\nu(\beta) = \nu(\alpha)+k \le n-1$ and 
  $\nu(\delta) \ge k$.  
  In this case, writing $\beta = ((1+p)^{r-1}-1)\beta'$, 
  we have
  \[ 
    x_r^{\alpha}z_r^{\beta} = x_r^{\alpha+(\varphi-1)\beta'}.
  \]
  Denote this element by $y$.  
  By Lemma \ref{symbolsr1}, we have that 
  $\ord(y,v) = 0$ and $\ord(y,u_0) \le n-k$.  
  Lemma \ref{unnondegen} forces $\ord(y,u_0) = n-k$, and so
  $f_n(y) = p^j$ by Proposition \ref{galcondr1}.
  The case statement now follows easily.
\end{proof}

\noindent
\footnotesize Dept.\ of Mathematics, Harvard University, 
Cambridge, MA\ \ 02138.\\
e-mail address: sharifi@math.harvard.edu
\end{document}